\documentclass[12pt]{article} 
\input amssym.def
\input amssym
\newtheorem{theorem}{Theorem}[section]

\newtheorem{algorithm}[theorem]{Algorithm}

\begin{document}
\def\edit{\bf }
\def\endedit{\rm }
\def\square{ }
\def\D{\Delta}
\def\L{\Lambda}
\def\Om{\Omega}
\def\S{\Sigma}
\def\s{\sigma}
\def\r{\rho}
\def\l{\lambda}
\def\b{\beta}
\def\a{\alpha}
\def\G{\Gamma}
\def\g{\gamma}
\def\k{\kappa}
\def\ve{\varepsilon}
\def\v{\varepsilon}
\def\vp{\varphi}
\def\om{\omega}
\def\proof{\noindent {\bf Proof: \quad}} 
 
\def\ppd{{\rm ppd}\,}
\def\lppd{{\rm lppd}\,}
\def\bppd{{\rm bppd}\,}
\def\Ker{{\rm Ker}\,}
\def\sim{{\rm sim}\,}

\def\Aut{{\rm Aut}\,}
\def\Out{{\rm Out}\,}
\def\Sym{{\rm Sym}\,}
\def\AGL{{\rm AGL}\,}
\def\PSL{{\rm PSL}\,}
\def\PSU{{\rm PSU}\,}
\def\PGL{{\rm PGL}\,}
\def\GF{{\rm GF}\,}
\def\GL{{\rm GL}\,}
\def\SL{{\rm SL}\,}
\def\Sp{{\rm Sp}\,}
\def\GSp{{\rm GSp}\,}
\def\GO{{\rm GO}}
\def\GU{{\rm GU}\,}
\def\U{{\rm U}\,}
\def\SU{{\rm SU}\,}
\def\O{{\rm O}\,}
\def\O+{{\rm O}$^+$\,}
\def\O-{{\rm O}$^-$\,}
\def\Oe{{\rm O}$^\v$\,}

\def\C{\mathcal{C}}
\def\calS{\mathcal{S}}

\def\qed{\hfill \rule{2mm}{2mm}}

\def\aut#1{{\rm{ Aut}}(#1)} 
\def\diff{\mathbin{\mkern-1.5mu\setminus\mkern-1.5mu}}

\def\la{\langle}
\def\ra{\rangle}
\def\G{\Gamma}
\def\Sym{\hbox{Sym\,}}
\def\soc{\hbox{soc\,}}

\title{Primitive prime divisor elements in finite classical groups}
\author{ Cheryl E. Praeger\thanks{Department of Mathematics, University
of Western Australia, Nedlands, W. A. 6907, Australia. } }
\date{}
\maketitle

\bigskip

\noindent{\bf Abstract}
This is an essay about a certain family of elements in the general 
linear group $\GL(d,q)$ called primitive prime divisor elements,
or $\ppd$-elements. A classification of the subgroups of $\GL(d,q)$ 
which contain such elements is discussed, and the proportions of
$\ppd$-elements in $\GL(d,q)$ and the various classical groups 
are given. This study of $\ppd$-elements was motivated by their 
importance for the design and analysis of algorithms
for computing with matrix groups over finite fields. An algorithm 
for recognising classical matrix groups, in which $\ppd$-elements
play a central role is described.
 
\section{Introduction.}
\label{sec-intro}%

The central theme of this essay is the study of a special kind of element
of the general linear group $\GL(d,q)$ of nonsingular $d\times d$ matrices 
over a finite field $\GF(q)$ of order $q$. We define these elements, 
which we call {\it primitive prime divisor elements} or 
$\ppd$-{\it elements}, and 
give good estimates of the frequencies with which they occur in $\GL(d,q)$
and the various classical matrix groups. Further we describe a classification of
the subgroups of $\GL(d,q)$ which contain $\ppd$-elements, and explore their
role in the design and analysis of a randomised algorithm for recognising 
the classical matrix groups computationally.

Perhaps the best way to introduce these ideas, and to explain the 
reasons for investigating this particular set of research questions, 
may be to give a preliminary discussion of a generic recognition 
algorithm for matrix groups. We wish to determine whether a given subgroup
$G$ of $\GL(d,q)$ contains a certain subgroup $\Om$. We design the algorithm
to study properties of randomly selected elements from $G$ in such a way
that, if $G$ contains $\Om$ then with high probability we will gain
sufficient information from these elements to conclude with certainty
that $G$ does contain $\Om$. A skeleton outline of the algorithm could
be written as follows.

\begin{algorithm} \label{basicalg}%
To recognise whether a given subgroup of $\GL(d,q)$
contains a certain subgroup $\Om$.
\begin{description}
\item[Input:] $G=\la X\ra\le\GL(d,q)$\quad and possibly some 
extra information about $G$.
\item[Output:] Either
\begin{description}
\item[(a)] $G$ contains the subgroup $\Om$, or
\item[(b)] $G$ does not contain $\Om$.
\end{description}
\end{description}
\end{algorithm}

If Algorithm~\ref{basicalg} returns option (a) then  $G$ 
definitely contains $\Om$. However if option (b) is returned there is 
a possibility that this response is incorrect. In other words  
Algorithm~\ref{basicalg} is a {\it Monte Carlo algorithm}.
It proceeds by making a sequence of random selections of 
elements from the group $G$, seeking a certain kind of subset $E$ of $G$,
which if found will greatly assist in deciding whether or not $G$ contains
$\Om$. The essential requirements for $E$ are two-fold:

\begin{description}
\item[1.]  If $G$ contains a subset $E$ with the required properties,
then either $G$ contains $\Om$, or $G$ belongs to a short list of
other possible subgroups of $\GL(d,q)$ (and the algorithm must then
distinguish subgroups in this list from subgroups containing $\Om$).
\item[2.] If $G$ contains $\Om$, then the event of {\it not} finding
a suitable subset $E$ in $G$ after a reasonable number $N(\v)$ of 
independent random selections of elements from $G$ has probability 
less than some small pre-assigned number $\v$.
\end{description}

In order to make the first requirement explicit, we need a classification 
of the subgroups of $\GL(d,q)$ which contain a
suitable subset $E$. 
Similarly in order to make the second requirement explicit, 
we need good estimates for the proportions of ``$E$-type elements'' 
in groups containing $\Om$.
Moreover, if these two requirements are to lead 
to an efficient algorithm for recognising whether $G$ contains $\Om$, the proportions of $E$-type elements in groups containing $\Om$ must be 
fairly large to guarantee that we have a good chance of finding
a suitable subset $E$ after a reasonable number of random selections;
and in practice we need good heuristics for producing approximately
random elements from a group.
Also, among other things, we need efficient procedures to identify
$E$-type elements, and to distinguish between the subgroups on the 
short list and the 
subgroups which contain $\Om$. The aim of this paper is to present 
and discuss results of these types, and the corresponding 
recognition algorithms, in the cases where
$\Om$ is one of the classical matrix groups. 
In these cases the subset $E$ consists of certain $\ppd$-elements.

I am grateful to Igor Shparlinski for some very helpful discussions
and advice on the analysis in Section~\ref{sec-complexity}.
Eamonn O'Brien made a careful reading of an early draft of the paper
and the current version has been much improved as a result of his 
detailed comments. Also I thank John Cannon for making available to 
me the results mentioned in Section~\ref{sec-performance} of some 
tests of the {\sc Magma} implementation of the classical recognition algorithm.

\section{Classical groups.}\label{sec-classical}%

We consider certain subgroups of $\GL(d,q)$ where
 $d$ is a positive integer and $q=p^a$, a power 
of a prime $p$, and we let $V$ denote the underlying vector space of 
$d$-dimensional row vectors 
over $\GF(q)$ on which $\GL(d,q)$ acts naturally. 

The classical groups
 preserve certain bilinear, sesquilinear or quadratic forms
on $V$. To describe them we adapt some notation from the book of Kleidman and 
Liebeck~\cite{kl}. A subgroup $G$ of $\GL(d,q)$ is said to {\it preserve 
a form $\k$ modulo scalars} if there exists a homo\-morphism $\mu:G\rightarrow
\GF(q)^\#$ such that, in the case of a bilinear or sesquilinear form,
$\k(ug,vg)=\mu(g)\cdot\k(u,v)$, or, in the case of a quadratic form, $\k(vg)=\mu(g)\cdot\k(v)$, for all $u,v\in V$ and $g\in G$. A matrix 
$g$ in such a group is said to {\it preserve $\k$ modulo scalars}, and
if $\mu(g)=1$ then $g$ is said to {\it preserve} $\k$.
We denote by $\D$ or $\D(V,\k)$ the group of all matrices in $\GL(d,q)$ 
which preserve $\k$ modulo scalars, and by $S$ the subgroup of
$\D$ consisting of those matrices which preserve $\k$ and which 
have determinant 1. 

The subgroup $\Om$ which we shall seek to recognise
is equal to $S$ unless $\k$ 
is a non-degenerate quadratic form, and in this latter case $\Om$ has index 
2 in $S$ and is the unique such subgroup of $S$. 
There are four families of subgroups which we shall consider, and
by a {\it classical group} in $\GL(d,q)$ we shall mean a subgroup $G$ which
satisfies $\Om\le G\le \D$, for $\Om,\ \D$ in one of these families. The
four families are as follows.

\begin{description}
\item[(i)] {\it Linear groups:} $\k=0,\ \D =\GL(d,q)$ and $\Om = 
\SL(d,q)$;
\item[(ii)] {\it Symplectic groups:} $d$ is even, $\k$ is a non-degenerate alternating bilinear form on $V$, $\D =\GSp(d,q)$ and $\Om 
=\Sp(d,q)$;
\item[(iii)] {\it Orthogonal groups:} $\k$ is a non-degenerate quadratic 
form on $V$, $\D = \GO^\v(d,q)$, and $\Om =\Om^\v(d,q)$, where 
$\v=\pm$ if $d$ is even, and $\v=\circ$ if $d$ is odd. If $d$ is odd then
also $q$ is odd since $\k$ is non-degenerate;
\item[(iv)] {\it Unitary groups:} $q$ is a square, $\k$ is a 
non-degenerate unitary form on $V$, that is a non-degenerate 
sesquilinear form with respect to the automorphism of $\GF(q)$ of order
2, $\D = \GU(d,q)$ and $\Om = \SU(d,q)$.
\end{description}

The books~\cite{kl, det} are good references for information about
the finite classical groups. 
 
\section{Primitive prime divisors and $\ppd$-elements.}\label{sec-ppds} 

Let $b, e$ be positive integers with $b>1$. 
A prime $r$ dividing $b^e -1$ is said to be a {\it primitive prime 
divisor} of $b^e - 1$ if $r$ does not divide $b^i -1$ for any $i$ such
that $1 \leq i < e$. It was proved by Zsigmondy~\cite{zsig} in 1892 that  
$b^e - 1$ has a primitive prime divisor unless either the pair $(b, e)$ 
is $(2, 6)$, or $e = 2$ and $b+1$ is a power of 2.
Observe that
\[
|\GL(d,q)|=q^{d\choose 2} \prod_{1\le i\le d} (q^i-1).
\]
This means that primitive prime divisors of $q^e-1$ for various 
values of $e\le d$ divide $|\GL(d,q)|$, and indeed divide $|\Om|$ for
various of the classical groups $\Om$ in $\GL(d,q)$. We define 
{\it primitive prime divisor elements}, sometimes
called $\ppd$-{\it elements}, in $\GL(d,q)$ to be those elements with
order a multiple of some such primitive prime divisor.
Thus we define an element $g\in\GL(d,q)$ to be a $\ppd(d,q;e)$-{\it element}
if its order $o(g)$ is divisible by some primitive prime divisor
of $q^e-1$.

Our interest is mainly in $\ppd(d,q;e)$-elements with $e>d/2$ and we shall 
describe in Section~\ref{sec-ppdsgps} a classification 
by Guralnick, Penttila, Saxl and the author 
in~\cite{ppds} of all subgroups of $\GL(d,q)$ containing such an element.
We shall henceforth reserve the term $\ppd$-elements to refer to elements
of $\GL(d,q)$ which are $\ppd(d,q;e)$-elements for some $e>d/2$.
Note that, if $g\in\GL(d,q)$ is a $\ppd(d,q;e)$-element with $e>d/2$, 
then there is a unique $g$-invariant $e$-dimensional subspace of the
underlying vector space $V$ on which $g$ acts irreducibly, 
and also the characteristic polynomial
for $g$ has an irreducible factor over $\GF(q)$ of degree $e$. While
neither of these two conditions is sufficient to guarantee that an element is a
$\ppd(d,q;e)$-element, it turns out that most elements satisfying
either of them are in fact $\ppd(d,q;e)$-elements.
In addition, a large proportion of elements in any of the classical
groups are $\ppd$-elements, and this fact has proved to be very important 
for the development of recognition algorithms for classical groups.

In 1974 Hering~\cite{her1} investigated subgroups of $\GL(d,q)$ containing $\ppd(d,q;d)$-elements. Such subgroups act irreducibly on $V$. 
Hering was interested in applications of these
results to geometry, in particular for constructing finite translation
planes. He was also interested in the link between such groups and 
finite affine 2-transitive permutation groups. 
If $G$ is a finite affine 2-transitive permutation group  
acting on a set $X$, then
$X$ may be taken as the set of vectors of a finite vector space,
say $V=V(d,q)$ of dimension $d$ over $\GF(q)$, and 
$G= N G_o$ where $N$ is the group of translations of $V$ and $G_o$
is a subgroup of $\GL(d,q)$ acting transitively on $V^\#$,
that is $G_o$ is a {\it transitive linear group}. 
Conversely if $G_o$ is a transitive linear group on $V$, and 
$N$ is the group of translations of $V$, then $ N G_o$ is
a 2-transitive permutation group of affine type acting on $V$.
Thus the problems of classifying finite affine 2-transitive groups,
and classifying finite transitive linear groups are equivalent.
Moreover if $G_o$ is transitive on $V^\#$ then $q^d-1$ divides $|G_o|$
so that $G_o$ contains a $\ppd(d,q;d)$-element. Hering's work led to
a classification of finite affine 2-transitive permutation groups, 
see~\cite{her2} and also~\cite[Appendix]{lieb}. In common with most of the 
classifications we shall mention related to $\ppd$-elements, 
this classification depends on the classification of the finite simple groups.
Merkt~\cite{merkt} extended Hering's work obtaining a better
description of certain of the subgroups of $\GL(d,q)$ containing a  
$\ppd(d,q;d)$-element. 

Dempwolff~\cite{demp} in 1987 began an investigation of subgroups of 
$\GL(d,q)$ containing a $\ppd(d,q;e)$-element for some $e\ge d/2$. 
His analysis is independent of the work of Aschbacher which we shall 
describe in the next section, and he made significant progress
on describing what we shall call (and shall define in the next section) 
the ``geometric subgroups'' containing such $\ppd$-elements. He also
did some work on the nearly simple examples. 
The classification in~\cite{ppds} of all subgroups of $\GL(d,q)$ 
containing a  $\ppd(d,q;e)$-element for some $e>d/2$ uses the work of
Aschbacher to guide both the analysis and the presentation of the
examples. 
Similar results may be obtained if the condition $e>d/2$ is 
relaxed, but their proofs become more technical.

\section{Aschbacher's classification of finite linear groups.}
\label{sec-asch}%

Aschbacher's description~\cite{asch} of subgroups of $\GL(d,q)$,
where $q=p^a$ with $p$ prime, has
been very influential both on the way problems concerning
linear groups are analysed and on the way results about such
groups are presented. Aschbacher
defined eight families of subgroups $\C_1,\ldots,\C_8$ of $\GL(d,q)$
as follows. These families are usually defined in terms of some
geometrical property associated with the action on the underlying 
vector space $V$, and in all cases maximal subgroups of $\GL(d,q)$
in the family 
can be identified. Subgroups of $\GL(d,q)$ in these families
are therefore called {\it geometric subgroups}.
We indicate in parentheses the rough structure of
a typical maximal subgroup in the family. 
Note that $Z$ denotes the subgroup of scalar matrices in $\GL(d,q)$.
Also, as in~\cite{kl}, we denote by $b$ a cyclic group of order $b$, 
and for a prime $r$ we denote by $r^{1+2c}$ an extraspecial group 
of that order.

\begin{description}
\item[$\C_1$] These subgroups act reducibly on $V$, and maximal 
subgroups in the family are the stabilisers of proper 
subspaces\ (maximal parabolic subgroups).
\item[$\C_2$] These subgroups act irreducibly but imprimitively 
on $V$, and maximal subgroups in the family are the stabilisers 
of direct sum decompositions $V=\oplus_{i=1}^t V_i$ with $\dim V_i = d/t$\
(wreath products $\GL(d/t,q)\wr S_t$).
\item[$\C_3$] These subgroups preserve on $V$ the structure
of a vector space over an extension field of $\GF(q)$, and 
maximal subgroups in the family are the stabilisers 
of extension fields of $\GF(q)$ of degree $b$, where $b$ is a prime
dividing $d$\ (the groups $\GL(d/b,q^b).b$).
\item[$\C_4$] These subgroups preserve on $V$ the structure of a
tensor product of subspaces, and maximal subgroups in the family 
are the stabilisers of decompositions $V= V_1\otimes V_2$\  
(central products $\GL(b,q)\circ \GL(c,q)$ where $d=bc$).
\item[$\C_5$] These subgroups preserve on $V$ the structure
of a vector space over a proper subfield of $\GF(q)$; such a
subgroup is said to {\it be realisable over a proper subfield}.  
The maximal subgroups in the family are the stabilisers modulo scalars
of subfields of $\GF(q)$ of prime index $b$ dividing $a$\ 
(central products $\GL(d,q^{1/b})\circ Z$).
\item[$\C_6$] These subgroups have as a normal subgroup an
$r$-group $R$ of symplectic type ($r$ prime) which acts absolutely
irreducibly on $V$, and maximal subgroups in the family 
are the normalisers of these subgroups,   
$(Z_{q-1}\circ R).\Sp(2c,r)$, where $d=r^c$ and $R$ is an 
extraspecial group $r^{1+2c}$, or if $r=2$ then $R$ may alternatively 
be a central product $4\circ 2^{1+2c}$.
\item[$\C_7$] These subgroups preserve on $V$ a tensor decomposition 
$V=\otimes_{i=1}^t V_i$ with $\dim V_i=c$, and maximal subgroups 
in the family are the stabilisers of such decompositions\   
($(\GL(c,q)\circ\ldots \circ\GL(c,q)).S_t$, where $d=c^t$).
\item[$\C_8$] These subgroups preserve modulo scalars a non-degenerate 
alternating, or sesquilinear, or quadratic form on $V$, and maximal 
subgroups in the family are the classical groups.
\end{description}

The main result of Aschbacher's paper~\cite{asch} (or 
see~\cite[Theorem~1.2.1]{kl}) 
 states that, for a subgroup $G$ of $\GL(d,q)$ which does
not contain $\SL(d,q)$,  either $G$ is a
geometric subgroup, or the socle $S$ 
of $G/(G\cap Z)$ is a nonabelian simple group, and the 
preimage of $S$ in $G$ is absolutely irreducible on $V$,
is not realisable over a proper subfield, and is not a classical
subgroup (as defined in Section~\ref{sec-classical}). The
family of such subgroups is denoted $\calS$, and subgroups
in this family will often be referred to as {\it nearly simple}
subgroups. Aschbacher~\cite{asch}
also defined families of subgroups of each of the classical 
subgroups $\D$ in $\GL(d,q)$, analogous to $\C_1,\ldots,\C_8,\calS$,
and proved that each subgroup of a classical group $\D$ which 
does not contain $\Om$ belongs to one of these families. 

\section{Linear groups containing $\ppd$-elements.}
\label{sec-ppdsgps}%

The analysis in~\cite{ppds} to determine the subgroups of 
$\GL(d,q)$ which contain a $\ppd(d,q;e)$-element for some 
$e>d/2$, was patterned on a similar analysis
carried out in~\cite{recog} to classify subgroups of $\GL(d,q)$
which contain both a $\ppd(d,q;d)$-element and a 
$\ppd(d,q;d-1)$-element. Moreover the results in~\cite{ppds} 
seek to give information
about the smallest subfield over which such a subgroup $G$ is
realisable modulo scalars. We say that $G$ is {\it realisable
modulo scalars} over a subfield $\GF(q_0)$ of $\GF(q)$ if $G$
is conjugate to a subgroup of $\GL(d,q_0)\circ Z$.

Suppose that $G\le\GL(d,q)$ and that $G$ contains a 
$\ppd(d,q;e)$-element for some $e>d/2$, and let $r$
be a primitive prime divisor of $q^e-1$ which divides $|G|$.
Suppose moreover that
$\GF(q_0)$ is the smallest subfield of $\GF(q)$ such that
$G$ is realisable modulo scalars over $\GF(q_0)$. 

There is a recursive aspect to the description in~\cite{ppds} 
of such subgroups $G$ which are geometric subgroups. 
For example, the reducible subgroups $G$ leave invariant some 
subspace or quotient space $U$ of $V$ of dimension $m\ge e$,
and the subgroup $G^U$ of $\GL(m,q)$ induced by $G$ in its
action on $U$ contains a $\ppd(m,q;e)$-element. 
In~\cite{ppds} no further
description is given of these examples, though extra information
may be obtained about the group $G^U$ by applying the
results recursively.

Although the classification of the geometric examples is
not difficult, care needs to be taken in order not to miss
some of them. 
For example, while at first sight it might appear that
a maximal imprimitive subgroup $\GL(d/t,q)\wr S_t$ (where $t>1$)
cannot contain a $\ppd(d,q;e)$-element since $r$ does
not divide $|\GL(d/t,q)|$, it is possible sometimes for
$r$ to divide $|S_t|=t!$, so that we
do have some examples in the family $\C_2$.

To understand how this can happen, observe that the defining 
condition for $r$ to be a primitive prime divisor of $q^e-1$, namely
that $e$ is the least positive 
integer $i$ such that $r$ divides $q^i-1$, is equivalent to the
condition that $q$ has order $e$ modulo the prime $r$.
Thus $r=ke+1\ge e+1$ for some $k\ge 1$.
Sometimes we can have $r=e+1$ (which satisfies $d/2<r\le d$)
and hence in these cases an imprimitive subgroup 
$\GL(1,q)\wr S_d$ will contain $\ppd(d,q;e)$-elements.

Both of the above observations come into play in describing the
examples in the family $\C_3$. Here either the prime $r=e+1=d$ and 
the group $G$ is conjugate to a subgroup of $\GL(1,q^d).d$,
or $e$ is a multiple of a prime $b$ where $b$ is a proper 
divisor of $d$ and, replacing $G$ by a conjugate if necessary, 
$G\le\GL(d/b,q^b).b$ such that $G\cap\GL(d/b,q^b)$ contains a
$\ppd(d/b,q^b;e/b)$-element.

After determination of the geometric examples there remains the 
problem of finding the nearly simple examples. So suppose 
that $G$ is nearly simple and $S\le G/(Z\cap G) \le\Aut S$ 
for some nonabelian simple group $S$. What we need is
a list of all possible groups $G$ together with the values of
$d, e$ and $q_0$. Although there is no
classification of all the nearly simple subgroups of $\GL(d,q)$
in general, it is possible to classify those which contain
a $\ppd(d,q;e)$-element. The reason we can do this is that,
for each simple group $S$, 
the presence of a $\ppd(d,q;e)$-element in $G$ leads to both
upper and lower bounds for $d$ in terms of the parameters of 
$S$ strong enough to lead to a complete classification. 

On the one hand $d$ is at least the minimum degree of a faithful 
projective representation of $S$ over a field of characteristic $p$,
and lower bounds are available for this in terms of the parameters 
of $S$. 
On the other hand we have seen that $r=ke+1\ge e+1\ge (d+3)/2$, and 
in all cases we may deduce that $r$ divides $|S|$. 
Moreover we have an upper bound on
the size of prime divisors of $S$ in terms of the parameters of $S$.
For some simple groups $S$ the upper and lower bounds for $d$ obtained
in this way conflict and we have a proof that there are no
examples involving $S$. In many cases however this line of argument
simply narrows down the range of possible values for $d$, $e$ and $r$.
Often there are examples involving $S$, but,
in order to complete the classification, we need to have 
more information about small dimensional representations of 
$S$ in characteristic $p$ than simply the lower bound for the dimension
of such representations.

For example if $S=A_n$ with $n\ge 9$ then 
 $d\ge n-2$ if $p$ divides $n$ and $d\ge n-1$ otherwise,
by~\cite{wag1, wag2, wag3}. Moreover $r\le n$, so $(d+3)/2\le n$,
and we obtain $n-2\le d\le 2n-3$ and
$r=e+1$. The upper bound for $d$ cannot 
be improved since we may have $r=n=e+1$ infinitely often.
Thus we need more information
about small dimensional representations of $A_n$ in characteristic $p$.
For $n\ge 15$ this is available from a combination of
results of James~\cite{james} and Wagner~\cite{wag3}. We see that 
the representations of $A_n$ and $S_n$ of dimension $n-1$ or $n-2$ 
are those coming from the deleted permutation module in the natural representation. These give an infinite family of examples with $q_0=p$.
All other faithful projective representations of $A_n$ have
dimension greater than the upper bound on $d$.
For the remaining cases, where $n<15$, special arguments are required,
making full use of information in~\cite{atlas,modat}.
The result of this analysis is an explicit list of examples 
for alternating groups $S$.

The list of examples of linear groups containing $\ppd$-elements 
can be found in~\cite[Section~2]{ppds} and is not reproduced here.
Note that completing the classification of the nearly simple examples for
classical groups $S$ over fields of characteristic different from $p$
involved proving new results about small dimensional representations
of such groups over fields of characteristic $p$.

\section{Various applications of the ``$\ppd$ classification''.}
\label{sec-applications}%

The classification of subgroups of $\GL(d,q)$ containing 
$\ppd$-elements has already been used in a variety of 
applications concerning finite classical groups. 
In particular the papers~\cite{gk, gs} make use of it
to answer questions concerning the generation of finite
classical groups, while in~\cite{lieb2} it is used to
show that the finite classical groups are characterised by
their orbit lengths on vectors in their natural modules.
Information about the invariant generation
of classical simple groups (see~\cite{recog2, shalev})
can be deduced from the classification (in~\cite{recog2},
or see Section~\ref{sec-basic-ideas})
of subgroups of classical groups containing two different 
$\ppd$-elements. (Elements $x_1,\ldots,x_s$ of a group $G$
are said to generate $G$ invariably if $\la x_1^{g_1},\ldots,
x_s^{g_s}\ra$ is equal to $G$ for all $g_1,\ldots,g_s\in G$.)

Similarly in~\cite{br} the $\ppd$ classification, or 
more accurately the more specialised classification based on 
it (and described in Section~\ref{sec-basic-ideas}), can be used 
to deal with the finite classical groups in an analysis of finite 
groups with the permutizer property. A group $G$
is said to have the {\it permutizer property} if, for every
proper subgroup $H$ of $G$, there is an element $g\in G\setminus H$ 
such that $H$ permutes with $\la g\ra$, that is $\la g\ra H=
H\la g\ra$. The main result of~\cite{br} is that all finite
groups with the permutizer property are soluble. The proof
consists of an examination of a minimal counterexample to
this assertion, and the $\ppd$ classification can be used to
show that the minimal counterexample cannot be an almost
simple classical group.

\section{Two different $\ppd$-elements in linear groups}
\label{sec-basic-ideas}

The principal application up to now of the classification
of linear groups containing $\ppd$-elements has been the
development by Niemeyer and the author in~\cite{recog2} 
of a recognition algorithm
for finite classical groups in their natural representation. 
The basic idea of this algorithm is as described in 
Section~\ref{sec-intro}. Given a subgroup $G$ of a classical
group $\D$ in $\GL(d,q)$ (as described in Section~\ref{sec-classical}), 
we wish to determine if $G$ contains the corresponding classical group
$\Om$. We do this by examining randomly selected elements
from $G$.
The elements of $G$ which we seek by random selection 
are $\ppd(d,q;e)$-elements for various values of $e>d/2$, 
and an appropriate set of such elements will form the subset $E$
mentioned in Section~\ref{sec-intro}. 

It turns out that the proportion of $\ppd(d,q;e)$-elements in
any of the classical groups is very high (as
shown in Section~\ref{sec-probs}),
so we are very likely to find such an element after a few
independent random selections from any subgroup of $\D$ 
which contains $\Om$. Suppose then that we have indeed found a
$\ppd(d,q;e)$-element in our group $G$, for some $e>d/2$. 
The $\ppd$-classification just described 
then provides a restricted list of possibilities for the
group $G$. The task is to distinguish subgroups containing
$\Om$ from the other possibilities, and this task is a
nontrivial one.

For the purposes of presenting the basic strategy, we 
assume that $G$ is irreducible on $V$ and that we have 
complete information about any $G$-invariant bilinear, 
sesquilinear or quadratic forms on $V$.
There are standard tests in practice which may be used to
determine whether $G$ is irreducible on $V$ and to
find all $G$-invariant forms (see~\cite{hr,parker}).
Note that in an implementation of the algorithm in~\cite{recog2}
a different protocol may be followed for deciding the stage 
at which to obtain this precise information about $G$.
 Nevertheless, we may and shall assume that $G$ 
does not lie in the Aschbacher classes $\C_1$ or $\C_8$. 
Then, having found a $\ppd(d,q;e)$-element in $G$ for some $e>d/2$,
the $\ppd$-classification would still allow the possibility that
$G$ lies in one of $\C_2, \C_3, \C_5, \C_6$, or that $G$ is 
nearly simple, as well as the desired conclusion that $G$ 
contains $\Om$. In the nearly simple case, the classification 
in~\cite{ppds} shows that there are approximately
30 infinite families and 60 individual examples of nearly
simple groups in explicitly known representations.

Guided by the original $\SL$-recognition algorithm developed
in~\cite{recog}, we decided to seek, in the first instance, 
{\it two different $\ppd$-elements} in $G$ by which we mean a
$\ppd(d,q;e)$-element and a $\ppd(d,q;e')$-element, where
$d/2<e<e'\le d$. We also decided to strengthen the $\ppd$-property
required of these elements in two different ways, by requiring
at least one of the $\ppd$-elements to be large and at least 
one of them to be basic.

Let $q=p^a$, and let $r$ be a primitive prime divisor of $q^e-1$.
Recall that $r=ke+1$ for some integer $k$. We say that $r$ is a {\it basic
primitive prime divisor} if $r$ is a primitive prime divisor of
$p^{(ae)}-1$, and that $r$ is a {\it large primitive prime divisor} 
if either $r\ge 2e+1$, or $r=e+1$ and $(e+1)^2$ divides $q^e-1$. 
Correspondingly we say that a $\ppd(d,q;e)$-element $g$ is {\it 
basic} if $o(g)$ is divisible by a basic primitive prime divisor 
of $q^e-1$, and that $g$ is {\it large} if $o(g)$ is divisible by  
a large primitive prime divisor $r$ of $q^e-1$ and either $r\ge 2e+1$ 
or $r=e+1$ and $(e+1)^2$ divides $o(g)$. 
Note that, for $e\ge 2$, if $q^e-1$ has a primitive prime divisor,
then $q^e-1$ has a basic primitive prime divisor unless $(q,e)=(4,3)$ 
or $(8,2)$. Similarly an explicit list can be given for pairs 
$(q,e)$ for which $q^e-1$ has a primitive prime divisor but does not
have a large primitive prime divisor (see~\cite{feit,her1} 
or~\cite[Theorem~2.2]{recog2}). Thus in most cases $q^e-1$ has 
both a large primitive prime divisor and a basic primitive prime divisor;
and many $\ppd$-elements will be both large and basic.

We shall see in  Section~\ref{sec-probs} that requiring the additional 
condition of being
large or basic does not alter significantly the very good upper 
and lower bounds we can give for the proportion of $\ppd$-elements
in subgroups of $\D$ containing $\Om$. 

Suppose that we now have $G\subseteq\D$ for 
some classical group $\D$ in $\GL(d,q)$, with $G$ irreducible on the
underlying vector space $V$, and suppose also that we have complete 
information about $G$-invariant forms so that we can guarantee that
$G$ is not contained in 
the class $\C_8$ of subgroups of $\D$. Further we suppose that
$G$ contains two different $\ppd$-elements, say a $\ppd(d,q;e)$-element 
$g$ and a $\ppd(d,q;e')$-element $h$, where $d/2<e<e'\le d$.

In~\cite[Theorem~4.7]{recog2}, Niemeyer and the author 
refined the classification 
in~\cite{ppds} to find all possibilities for the group $G$. 
These possibilities comprise groups containing $\Om$, members
of the Aschbacher families $\C_2, \C_3$ and $\C_5$, and 
some nearly simple examples. The presence of two different 
$\ppd$-elements certainly restricts the possibilities within these 
families, but it is still difficult 
to distinguish some of them from groups containing $\Om$.

If we require that at least one of $g, h$ is large and at least 
one is basic then, as was shown in~\cite[Theorem~4.8]{recog2}, 
the possibilities for irreducible subgroups
$G$ which do not contain $\Om$ are certain
subgroups in $\C_3$ and nearly simple groups in a very short list
comprising explicit representations of one infinite family and 
five individual nearly simple groups.

After our discussion of the proportions of $\ppd$-elements in 
classical groups in Section~\ref{sec-probs},  we shall return to
our discussion of the recognition algorithm. We shall see that the
algorithm can be completed by simply seeking a few more
$\ppd$-elements of a special kind which, if found, will rule out all 
but one possibility for $G$, enabling us to conclude that $G$ 
contains $\Om$.
 
\section{Proportion of $\ppd$-elements in classical groups.}
\label{sec-probs}

The questions we wish to answer from our discussion in this section are
the following. If $\Om\le G\le\D\le\GL(d,q)$, and $G$ contains two
different $\ppd$-elements at least one of which is large and at least
one of which is basic, then what is the probability of finding two
such elements after a given number $N$ of independent random selections
of elements from $G$? In particular, for a given positive real
number $\v$, is it
true that the probability of failing to find such elements after $N$ 
selections is less than $\v$ provided $N$ is sufficiently large?
And if so just how large must $N$ be?

These questions can be answered using simple probability theory
provided that we can determine, for a given $e$ (where $d/2<e\le d$), 
the proportion $\ppd(G,e)$ of elements of $G$ which are 
$\ppd(d,q;e)$-elements. This proportion may depend on the nature of 
the classical group
$\D$: that is, on whether $\D$ is a linear, symplectic, orthogonal 
or unitary group. 
In particular $\ppd(G,e)=0$ if $\D$ is a symplectic or orthogonal 
group and $e$ is odd, or if $\D$ is a unitary group and $e$ is even,
or if $\D$ is of type ${\rm O}^+$ and $e=d$.
This can be seen easily by examination of the orders of these groups. In all
other cases, provided that $d$ and $q$ are not too small, any subgroup
of $\D$ which contains $\Om$ will contain $\ppd(d,q;e)$-elements.

So suppose now that $\Om\le G\le\D$, that $d/2<e\le d$, and that 
$G$ contains a $\ppd(d,q;e)$-element $g$. It is not difficult 
(see~\cite[Lemma~5.1]{recog2}) to show that $V$ has a unique 
$e$-dimensional $g$-invariant subspace $W$ and that $g$ acts 
irreducibly on $W$. Moreover, if $\D$ is a symplectic, orthogonal, 
or unitary group, then $W$ must be nonsingular with respect to the
bilinear, quadratic, or sesquilinear form defining $\D$.

Next (see~\cite[Lemma~5.2]{recog2}) we observe that the group 
$G$ acts transitively on the set of all
nonsingular $e$-dimensional subspaces of $V$ (or all $e$-dimensional
subspaces if $\D=\GL(d,q)$). Thus the proportion of 
$\ppd(d,q;e)$-elements in $G$ is the same as the proportion of 
such elements which fix a particular nonsingular $e$-dimensional 
subspace $W$. Therefore we need to determine the proportion of 
$\ppd(d,q;e)$-elements in the setwise stabiliser $G_W$ of $W$ in $G$.

Now consider the natural map $\vp: g\mapsto g|_W$ which sends $g\in G_W$
to the linear transformation of $W$ induced by $g$. Then $\Om(W)\le
\vp(G)\le\D(W)\le\GL(W)$, and $\D(W)$ has the same type (linear, symplectic,
orthogonal, or unitary) as $\D$. If $g\in G_W$ and $g$ is a
$\ppd(d,q;e)$-element, then every element of the coset $g\Ker\vp$ is
also a $\ppd(d,q;e)$-element, since all elements in the coset
induce the same linear transformation $g|_W$ of $W$. Moreover 
in this case $g|_W$ is a $\ppd(e,q;e)$-element in $\vp(G)$ and
all such elements arise as images under $\vp$ of 
$\ppd(d,q;e)$-elements in $G_W$. 
It follows that $\ppd(G,e)$ is equal to the proportion 
$\ppd(\vp(G),e)$ of $\ppd(e,q;e)$-elements in $\vp(G)$.

Thus it is sufficient for us to determine $\ppd(G,d)$ for each 
of the possibilities for $\D$ which contain $\ppd(d,q;d)$-elements. 
This was done already by Neumann and the author 
in~\cite[Lemmas~2.3 and 2.4]{recog} in the case where $\D=\GL(d,q)$.
The techniques used there work also in the other cases although
some care is needed. 
The basic ideas are as follows. 

Let $g$ be a $\ppd(d,q;d)$-element in $G$, and let $C:=C_G(g)$.
Then $C$ is a cyclic
group, called a Singer cycle for $G$, and has order $n$ say, where
$n$ divides $q^d-1$ and $n$ is divisible by some primitive prime
divisor of $q^d-1$. The group $C$ is self-centralising in $G$.
Further each $\ppd(d,q;d)$-element in $G$ lies in a unique $G$-conjugate
of $C$. The number of $G$-conjugates of $C$ is $|G:N_G(C)|$, and
so the number of $\ppd(d,q;d)$-elements in $G$ is equal to 
$|G:N_G(C)|$ times the number of such elements in $C$. It follows that 
$$
\ppd(G,d)=|G:N_G(C)|\cdot \ppd(C,d)\cdot {|C|\over |G|} = {\ppd(C,d)\over u},
$$
where $\ppd(C,d)$ is the proportion of $\ppd(d,q;d)$-elements in $C$,
and $u:=|N_G(C):C|$. In the linear, symplectic and unitary cases
$u=d$, while in the orthogonal case $u$ is either $d$ or $d/2$ depending
on which intermediate subgroup $G$ is ($\Om\le G\le\D$). In the orthogonal
case we certainly have $u=d$ if $G$ contains ${\rm O}\, (V)$.

Thus we need to estimate $\ppd(C,d)$. Let $\Phi$ denote the product of all
the primitive prime divisors of $q^d-1$ (including multiplicities), so that
$(q^d-1)/\Phi$ is not divisible by any primitive prime divisor of
$q^d-1$. In all cases $\Phi$ divides $n=|C|$. Moreover an element
$x\in C$ is a $\ppd(d,q;d)$-element if and only if $x^{n/\Phi}\ne 1$,
that is if and only if $x$ does not lie in the unique subgroup of
$C$ of order $n/\Phi$. Hence 
$$
\ppd(C,d)={n-n/\Phi\over n} = 1-{1\over\Phi},
$$
and therefore 
$$
\ppd(G,d)={1\over u}\ (1-{1\over\Phi}) < {1\over u}.
$$
Since each primitive prime divisor of $q^d-1$ is of the form
$kd+1\ge d+1$, the quantity $\Phi$ is at least $d+1$, and hence
$$
\ppd(G,d)\ge {1\over u}\ (1-{1\over d+1})
$$
so we have 
$$
{1\over u}\ ({d\over d+1})\ \le\ppd(G,d)< {1\over u}.
$$
Putting all of this together we see that in almost all cases 
$\ppd(G,d)$ lies between $1/(d+1)$ and $1/d$, with the exception
being some orthogonal cases where $\ppd(G,d)$ lies between
$2/(d+1)$ and $2/d$.

To pull back this result to the general case where $d/2<e\le d$, we
need to have some particular information about the group $\vp(G)$ in the 
orthogonal case in order to know which of the bounds apply.
It turns out (see~\cite[Theorem~5.7]{recog2})
that for all cases, and all $e$ for which $d/2<e\le d$ and 
$\D$ contains $\ppd(d,q;e)$-elements, we have
$$
{1\over e+1}\le\ppd(G,e)<{1\over e}
$$
except if $\D$ is an orthogonal group of minus type, $e=d$ is even, 
and $G\cap {\rm O}^-(d,q)$ is either $\Om^-(d,q)$ (for any $q$) or
${\rm SO}^-(d,q)$ (for $q$ odd), in which case $2/(d+1)\le\ppd(G,d)<
2/d$.

Further (see~\cite[Theorem~5.8]{recog2}), the proportion 
of large $\ppd(d,q;e)$-elements in $G$
and the proportion of basic $\ppd(d,q;e)$-elements in $G$, 
whenever such elements exist, also lie
between the lower and upper bounds we have above for $\ppd(G,e)$.

In the classical recognition algorithm in~\cite{recog2} we are not
especially interested at first in particular values of $e$. We simply
wish to find $\ppd$-elements for some $e$ between $d/2$ and $d$.
The proportion of such elements in $G$ is
$$
\ppd(G):=\sum_{d/2<e\le d} \ppd(G,e).
$$
In the linear case, where $\Delta =\GL(d,q)$, this is approximately equal to
$\sum_{d/2<e\le d} e^{-1}$ which, in turn, is approximately 
$$
\int_{d/2}^d{dx\over x} = \log 2 = 0.693\ldots
$$
while in the other cases $\ppd(G)$ is approximately equal to the sum
of $e^{-1}$ either over all even $e$, or all odd $e$ between $d/2$ 
and $d$; this is approximately equal to $(\log 2)/2$. These
computations can be done carefully resulting in very good upper and
lower bounds for $\ppd(G)$ which differ by a small multiple of $d^{-1}$,
see~\cite[Theorem~6.1]{recog2}. Moreover, except for small values of $d$,
these upper and lower bounds for $\ppd(G)$ are also upper and lower 
bounds for the proportions of large $\ppd$-elements and of basic 
$\ppd$-elements in $G$.

We can model the process of random selection of $N$ elements from $G$, 
seeking $\ppd$-elements, as a sequence of $N$ binomial trials with
probability of success on each trial (that is, each selection) being $\ppd(G)$. 
Using this model we can compute the probability of finding (at
least) ``two different $\ppd$-elements'' after $N$ 
independent random selections.
The extent to which this computed probability measures the true
probability in a practical implementation depends on whether the 
assumptions for the binomial model hold for the implementation. In
particular the binomial model will give a good fit if the selection
procedure is approximately {\it uniform}, that is the probability 
of selecting each element of $G$ on each selection is 
approximately $|G|^{-1}$, and if the selections are approximately 
independent. For any small positive real number $\ve$, under the 
binomial model the probability of failing to find ``two different
 $\ppd$-elements'' after $N$ independent uniform random selections
is less than $\ve$ provided that $N$ is greater than a small (specified)
multiple of $\log \ve^{-1}$, see~\cite[Theorem~6.4 and Lemma~6.5]{recog2}.

The same approach (under the same assumptions about uniformity 
and independence of the random selections) gives good estimates for the 
number $N=N(\v)$ of selections needed in order that  the probability 
of failing to find ``two different $\ppd$-elements'', at least 
one of which is large and at least one of which is basic, after 
$N$ random selections is less than $\ve$. Namely $N(\ve)$ is a small 
(specified) multiple of $\log \ve^{-1}$. For example, if 
$\Delta =\GL(d,q)$ with $40\le d\le 1000$ and $\ve= 0.1$, then
$N(\ve)=5$.

\section{Classical recognition algorithm: an outline}
\label{sec-outline}

Suppose that $G\subseteq\D$ for 
some classical group $\D$ in $\GL(d,q)$, with $G$ irreducible on the
underlying vector space $V$, and that we have complete 
information about $G$-invariant forms (so that $G$ is not contained in 
the class $\C_8$ of subgroups of $\D$). We wish to determine whether
or not $G$ contains the corresponding classical group $\Om$.
Our algorithm is a Monte Carlo algorithm which may occasionally
fail to detect
that $G$ contains $\Om$. The probability of this happening is
less than a predetermined small positive real number $\ve$.

First we make a number $N$ of independent uniform random
selections of elements from $G$, where $N\ge N(\ve/3)$ as in 
Section~\ref{sec-probs}. If we fail to find two different
$\ppd$-elements in $G$, with at least one of them large and at least 
one basic, then we report that $G$ does not contain
$\Om$. There is a possibility that this response is incorrect,
but if in this case $G$ does contain $\Om$ then from Section~\ref{sec-probs},
the probability of failing to find suitable elements is less than $\ve/3$.
Thus the probability of reporting at this stage that  $G$ does not contain
$\Om$, given that  $G$ does contain $\Om$, is less than $\ve/3$.

Suppose now that
$G$ contains two different $\ppd$-elements, say a $\ppd(d,q;e)$-element 
$g$ and a $\ppd(d,q;e')$-element $h$, where $d/2<e<e'\le d$, and that
at least one of $g, h$ is large and at least one is basic. As 
discussed in Section~\ref{sec-basic-ideas}, the possibilities for $G$ are
that (i)~$G\supseteq\Om$, or that (ii)~$G$ is conjugate to a subgroup of
$\GL(d/b,q^b).b$ for some prime $b$ dividing $d$, or that (iii)~$G$ is
one of a very restricted set of nearly simple groups. 
In order to distinguish case
(i) from cases (ii) and (iii) it turns out that essentially we need to
find a few extra $\ppd$-elements. 

The ``extension field groups'' in case (ii) are the most difficult to 
handle. The basic idea here can be illustrated by considering the 
linear case where $\Delta=\GL(d,q)$. For a prime $b$ dividing $d$, 
the only values of $e$ such that $\GL(d/b,q^b).b$ contains a 
$\ppd(d,q;e)$-element are multiples of $b$ (apart from the exceptional case
where $b=d$ and $d$ is a primitive prime divisor of $q^{d-1}-1$).
Thus finding in $G$ a $\ppd(d,q;e)$-element for some $e$ which is not
a multiple of $b$ will prove that $G$ is not conjugate to a subgroup 
of $\GL(d/b,q^b).b$. If $G\supseteq\Om$, then the proportion of such 
elements in $G$ is $\ppd(G) - \sum_{d/2<ib\le d}\ppd(G,ib)$ which is
approximately equal to $\ppd(G) - (\sum_{d/(2b)<i\le d/b} (ib)^{-1})$.
This in turn is approximately equal to $\log 2 - b^{-1}\log 2 = (\log 2)
(b-1)/b$. By~\cite[Theorem~8.30]{nzm}, the number $\mu(d)$ of 
distinct primes dividing $d$ is
$O(\log d/\log\!\log d)$. Arguing as in Section~\ref{sec-probs}, 
there is an integer $N_b(\ve)$ such that, if $G\supseteq\Om$, 
then the probability of failing to find a $\ppd(d,q;e)$-element in $G$
with $e$ coprime to $b$ after $N_b(\ve)$ independent random selections
is less than $\ve/3\mu(d)$. If $G\supseteq\Om$, then we may need to find up to
$\mu(d)$ extra $\ppd$-elements to eliminate case (ii) as a possibility,
and the probability of failing to eliminate it after $N$ random selections
from $G$, where $N$ is the maximum of the $N_b(\ve)$, is less than $\ve/3$.
If we fail to find the required set of elements after these $N$ further
random selections then we report that $G$ does not contain $\Om$.
Thus the probability of reporting at this second 
stage that  $G$ does not contain
$\Om$, given that  $G$ does contain $\Om$, is less than $\ve/3$.
The number $N$ of selections we need to make for this second stage
is $O(\log \ve^{-1} + \log\log d)$.
Eliminating possibility (ii) for the other classical groups
is done using these basic ideas, but the details are considerably
more complicated for the symplectic and orthogonal groups when $b=2$.

For each of the nearly simple groups which contain two different
$\ppd$-elements $g, h$ as above, there are in fact only two values of
$e$ for which the group contains $\ppd(d,q;e)$-elements, namely
the values corresponding to the elements $g$ and $h$. To distinguish
groups $G$ containing $\Om$ from this nearly simple group we simply 
need to find in $G$ a $\ppd(d,q;e)$-element for a third value of $e$.
For each pair $(d,q)$ there is only a small number of possible
nearly simple groups (usually at most 1, and in all cases at most 3).
As before there is some $N_{\sim}(\ve)$ such that, if $G\supseteq\Om$,
then the probability of failing to find suitable elements to eliminate
these nearly simple groups after $N_{\sim}(\ve)$ random selections
from $G$ is less than $\ve/3$.
If we fail to find the required elements after $N_{\sim}(\ve)$ further
random selections then we report that $G$ does not contain $\Om$.
Thus the probability of reporting at this third and final  
stage that  $G$ does not contain
$\Om$, given that  $G$ does contain $\Om$, is less than $\ve/3$.

Once we have found all the required elements to remove possibilities 
(ii) and (iii) we may report with certainty that
$G$ does contain $\Om$. 

The probability that the algorithm reports that $G$ does 
not contain $\Om$, given that  $G$ does contain $\Om$, is less than $\ve$.
The requirements to bound the probability of error at the 
three stages of the algorithm are such that the complete algorithm 
requires us to make
$O(\log\ve^{-1} + \log\!\log d)$ random selections from $G$.

\section{Computing with polynomials}
\label{sec-polynomials}

In this section we describe how we process an element $g$ of a classical group 
$\D\le\GL(d,q)$ to decide if it is a $\ppd$-element. This is a 
central part of the algorithm. 

The first step is to compute the characteristic polynomial $c_g(t)$ 
of $g$, and to determine whether or not $c_g(t)$ has an irreducible 
factor of degree greater than $d/2$. If no such factor exists then
$g$ is not a $\ppd$-element. So suppose that $c_g(t)$ has an 
irreducible factor $f(t)$ of degree $e>d/2$. 

Thus we know that there is a unique $g$-invariant $e$-dimensional
subspace $W$ of $V$ and that the linear transformation $g|_W$ induced
by $g$ on $W$ has order dividing $q^e-1$; $g$ will be a 
$\ppd(d,q;e)$-element if and only if the order of $g|_W$ is divisible by
some primitive prime divisor of $q^e-1$. By an argument introduced 
in Section~\ref{sec-probs}, this will be the case if and only if 
$(g|_W)^{(q^e-1)/\Phi}\ne 1$, where $\Phi = \Phi(e,q)$ and $\Phi(e,q)$ denotes 
the product of all the primitive prime divisors of $q^e-1$ 
(including multiplicities).
Determining whether or not this is the case can be achieved by
computing within the polynomial ring $\GF(q)[t]$ modulo the
ideal $\la f(t)\ra$, namely $(g|_W)^{(q^e-1)/\Phi}$ will be a
non-identity matrix if and only if $t^{(q^e-1)/\Phi}\ne 1$ in
this ring.

We can test whether of not $g$ is a large or basic
$\ppd(d,q;e)$-element by the same method with $\Phi(e,q)$ replaced by
$\Phi_l(e,q)$ or $\Phi_b(e,q)$ respectively. Here $\Phi_l(e,q)$
and $\Phi_b(e,q)$ are the products of all the large and basic 
primitive prime divisors of $q^e-1$ (including multiplicities) 
respectively. 

The idea for checking the $\ppd$-property by
determining whether a single power of $g$ is 
the identity comes from the special linear recognition algorithm 
in~\cite{recog}, while the idea of deciding this by a
computation in the polynomial ring is due to Celler and  
Leedham-Green~\cite{cl}.

\section{Complexity of the classical recognition algorithm}
\label{sec-complexity}

In~\cite[Section~4]{dimacs} an analysis of the running cost for
the classical recognition algorithm was given based on 
``classical'' algorithms for computing in finite fields.
For example the cost of multiplying two $d\times d$ matrices
 was taken to be $O(d^3)$ field operations (that is, additions,
multiplications, or computation of inverses). We take this
opportunity to re-analyse the algorithm in terms of more modern
methods for finite field computations. These methods can lead to
improvements in performance over the classical methods. However
efficient implementation of the modern methods is a highly
nontrival task requiring substantial effort, see for example
the paper of Shoup~\cite{sh} which addresses the problem of 
efficient factorisation of polynomials over finite fields.
I am grateful to Igor 
Shparlinski for some interesting and helpful discussions
concerning such algorithms. 

The {\it exponent of matrix multiplication} is defined as the infimum
of all real numbers $x$ for which there exists a 
matrix multiplication algorithm which requires no more than
$O(d^x)$ field operations to multiply together two 
$d\times d$ matrices over a field of order $q$. It
is denoted by $\om$ or $\om(d,q)$.
Thus, for all positive real numbers $\ve$, there exists
such an algorithm which requires $O(d^{\om +\ve})$ field operations, that
is matrix multiplication can be performed with $O(d^{\om+o(1)})$
field operations. 
In~\cite[Sections~15.3,~15.8]{bcs} an algorithm is given and analysed
for which $O(d^x)$ field operations are used with $x < 2.39$ (and hence $\om<2.39$),  and it was shown there also that
$\om$ can depend (if at all) only on the prime $p$ dividing $q$
rather than on the field size $q$.
Moreover the cost of performing a field operation depends on the 
data structure used to represent the field and is 
$O((\log q)^{1+o(1)})$ for each field operation, that is, the cost is 
$O((\log q)^{1+\ve})$ for each $\ve >0$. 

Now let $\mu$ be the cost of producing a single random element from
the given subgroup $G=\la X\ra$ of $\GL(d,q)$. As discussed
in~\cite[p.~190]{dim1}, theoretical methods for producing
approximately random elements from a matrix group are not good
enough to be translated into practical procedures for 
use with algorithms such as the classical recognition algorithm.
For example, Babai~\cite[Theorem~1.1 and Proposition~7.2]{babai} produces,
from a given generating set $X$ for a subgroup $G\le\GL(d,q)$, 
a set of $O(d^2\log q)$ elements of $G$ at a cost of $O(d^{10}
(\log q)^5)$ matrix multiplications, from which nearly uniformly 
distributed random elements of $G$ can be produced at a cost of
$O(d^{2}\log q)$ matrix multiplications per random element.
The practical implementation of the classical recognition algorithm
uses an algorithm developed in~\cite{clmno} for producing 
approximately random elements in classical groups which, when
tested on a range of linear and classical groups was found to produce,
for each relevant value of $e$, $\ppd(d,q;e)$-elements in 
proportions acceptably close to the true proportions in the group.
This procedure has an initial phase which costs $O(d^{\om+o(1)})$ field
operations, and then the cost of producing each random element is
$O(d^{\om+o(1)})$ field operations (see also~\cite[Section~4.1]{dimacs}).
Further analysis of the algorithm in~\cite{clmno} may be found 
in~\cite{cg, ds1, ds2}.

Testing each random element $g\in G$ 
involves first finding its characteristic polynomial 
$c_g(t)$. The cost of doing this deterministically is
$O(d^{\om+o(1)})$ field operations (see~\cite{k-g} or~\cite[Section~16.6]{bcs}).
Next we test whether $c_g(t)$ has an irreducible factor of degree
greater than $d/2$. This can be done deterministically at a cost
of $O(d^{\om+o(1)} + d^{1+o(1)}\log q)$ field operations, see~\cite{ks}.
(Although the full algorithm in~\cite{ks} for obtaining a complete 
factorisation of $c_g(t)$ is non-deterministic, we only need the first two
parts of the algorithm, the so-called square-free factorisation 
and distinct-degree factorisation procedures, and these are 
deterministic.) 
Suppose now that $c_g(t)$ has an irreducible factor $f(t)$ of 
degree $e>d/2$. We then need to compute $\Phi(e,q)$, the product of all
the primitive prime divisors of $q^e-1$ (counting multiplicities).
A procedure for doing this is given in~\cite[Section~6]{recog}.
It begins with setting $\Phi = q^e-1$ and proceeds by repeatedly dividing
$\Phi$ by certain integers. The procedure runs over all the
distinct prime divisors $c$ of $e$, and by~\cite[Theorem~8.30]{nzm}
there are $O(\log e/\log\log e)= O(\log d/\log\log d)$ such prime 
divisors. For each $c$, the algorithm computes twice the greatest 
common divisor of two positive integers where
the larger of the two integers may be as much as $q^e$, and then
makes up to $d\log q$ greatest common divisor computations 
for which the larger of the two integers
is $O(d)$. By~\cite[Theorem~8.20 and its Corollary]{ahu} (or see
\cite[Note 3.8]{bcs}), the cost of computing the greatest common 
divisor of two positive integers less than $2^n$, is $O(n(\log n)^{O(1)})$
bit operations. It follows that the cost of computing 
$\Phi(e,q)$ is $O(d(\log d)^{O(1)} (\log q)^2)$ bit operations.
Having found $\Phi(e,q)$, we need to determine whether  
$t^{(q^e-1)/\Phi(e,q)}$ is equal to $1$ in the polynomial ring $\GF(q)[t]$
modulo the ideal $\la f(t)\ra$. This involves  $O(d\log q)$ multiplications 
modulo $f(t)$ of two polynomials of degree less than $d$ over $\GF(q)$. Each of these polynomial multiplications costs $O(d\log d\log\log d)$
field multiplications, (see~\cite[Theorem~2.13 and Example~2.6]{bcs}). 
Thus this test costs 
$O(d^2\log d\log\log d\log q)$ field operations.
Therefore the cost of testing whether a random element $g$ is a $\ppd$-element
is $O(d^{\om+o(1)} + d^2\log d\log\log d\ \log q)$ field operations plus
$O(d (\log d)^{O(1)} (\log q)^2))$ bit operations, and hence is
$$O(
d^{\om+o(1)} (\log q)^{1+o(1)}+ 
d^2\log d\log\log d\, (\log q)^{2+o(1)} )
$$
bit operations. This is at most $O(d^{\om+o(1)} 
(\log q)^{2+o(1)})$ bit operations.
The cost of checking whether $g$ is a large $\ppd$-element
is the same as this. To check if $g$ is a basic $\ppd(d,q;e)$-element
involves computing $\Phi_b(e,q)=\Phi(ae,p)$ (where $q=p^a$) 
instead of $\Phi(e,q)$. Arguing as above, the cost of computing 
$\Phi_b(e,q)$ is $O(ad(\log (ad))^{O(1)} (\log p)^2) = 
O(d(\log d)^{O(1)} (\log q)^2)$
bit operations, and hence the cost of testing whether $g$ is a basic
$\ppd$-element is also at most $O(d^{\om+o(1)} (\log q)^{2+o(1)})$
bit operations. 

Since we need to test $O(\log \ve^{-1} + \log\log d)$ elements of $G$, the
total cost of the algorithm is as follows.

 \begin{theorem}
Suppose that $G\subseteq\D$ for 
some classical group $\D$ in $\GL(d,q)$, with $G$ irreducible on the
underlying vector space $V$, and that we have complete 
information about $G$-invariant forms (so that $G$ is not contained in 
the class $\C_8$ of subgroups of $\D$). Assume that $d$ is large 
enough that $\Om$ contains two different $\ppd$-elements with at 
least one of them large and at least one basic.
Further let $\ve$ be a positive real number with $0<\ve <1$.
Assume that we can make uniform independent random selections 
of elements from $G$ and that the cost of producing each random 
element is $\mu$ bit operations.
Then the classical recognition algorithm in~\cite{recog2} 
uses $O(\log \ve^{-1} + \log\log d)$ random elements from $G$ to test
whether $G$ contains $\Om$, and in the case where $G$ contains $\Om$,
the probability of failing to report that $G$ contains $\Om$ is less
than $\ve$. The cost of this algorithm is
$$
O\big( 
(\log \ve^{-1} + \log\log d) (\mu + d^{\om+o(1)} (\log q)^{2+o(1)})                      \big)
$$
bit operations, where $\om$ is the exponent of matrix multiplication. \end{theorem}

\section{Classical recognition algorithm: final comments}
\label{sec-performance}

The classical recognition algorithm in~\cite{recog2} 
has been implemented and is available as part of the 
{\it matrix} share package with the {\sf GAP} system~\cite{gap},
and is also implemented in {\sc Magma}~\cite{magma}.
In the {\sc Magma} implementation rather large groups 
have been handled by the algorithm without problems: 
John Cannon has informed us that, on a SUN Ultra 2 
workstation with a 200 MHz processor, recognising 
$\SL(5000, 2)$, for example, took 3214 CPU seconds averaged over
five runs, while recognising $\SL(10000, 2)$ was possible in 
14334 CPU seconds, again averaged over five runs of the algorithm.

The algorithm as described in this paper relies on the presence in 
the classical group $\Om$ of two different $\ppd$-elements, where 
at least one is large and at least one is basic. However, for some
small values of the dimension $d$, depending on the type of the 
classical group and the field order $q$, $\Om$ may not contain
such elements. In these cases a modification of the algorithm
has been produced in~\cite{recog3} which makes use of
elements which are similar to $\ppd$-elements. The results in
\cite{recog3} demonstrate that, with some effort, it is possible to
extend the probability computations in Section~\ref{sec-probs}.

An alternative algorithm to recognise classical groups in their 
natural representations has been developed by Celler and Leedham-Green
in~\cite{cl2}. This algorithm also uses the
Aschbacher classification~\cite{asch} of subgroups of $\GL(d,q)$ as
its organisational principle. Like the algorithm in~\cite{recog2} it 
makes use of a search by random selection for certain elements.
Although no analysis of the complexity of the algorithm is given
in~\cite{cl2}, the analysis we give in Section~\ref{sec-complexity}
gives a reasonable measure of the complexity of this algorithm also.
Finally, as with the algorithm in~\cite{recog2}, the algorithm 
in~\cite{cl2} does not work for certain families of small dimensional
classical groups (notably the groups of type ${\rm O}^+(8,q)$), and the
methods of~\cite{recog3} are required to deal with these groups.


\begin{thebibliography}{99}
\bibitem{ahu}%
A. V. Aho, J. E. Hopcroft and J. D. Ullman,
{\it The design and analysis of computer algorithms},
Addison-Wesley, 1975.

\bibitem{asch}%
M. Aschbacher, On the maximal subgroups of the finite 
classical groups, {\it Invent. Math.} {\bf 76} (1984), 469-514. 

\bibitem{babai}%
L. Babai,
Local expansion of vertex-transitive graphs and random generation
in finite groups, {\it Proc. $23$rd ACM STOC} (1991), 164--174.

\bibitem{br}%
J. C. Beidleman and D. J. S. Robinson, 
On finite groups satisfying the permutizer condition,
{\it J. Algebra} {\bf 191} (1997), 686--703.

\bibitem{magma}%
W. Bosma and J. J. Cannon, Handbook of
{\sc Magma} functions, {\it Department of Pure Mathe\-ma\-tics, Sydney
University}, 1993.

\bibitem{bcs}%
P. B\"urgisser, M. Clausen and M. A. Shokrollahi,
{\it Algebraic complexity theory}, Springer, Berlin Heidelberg, 1997.

\bibitem{cl}%
F. Celler and C.R. Leedham-Green,
Calculating the order of an invertible matrix, in  
{\it Groups and Computation II}, DIMACS: Series in Discrete 
Mathematics and Theoretical Computer Science, {\bf 28}  (ed. 
L. Finkelstein and W. M. Kantor, American Mathematical Society, 
1997) 55--60. 

\bibitem{cl2}%
F. Celler and C.R. Leedham-Green,
A non-constructive recognition algorithm for the special linear and other
classical groups, in  
{\it Groups and Computation II}, DIMACS: Series in Discrete Mathematics and
Theoretical Computer Science, {\bf 28}  (ed. L. Finkelstein and W. M.
Kantor, American Mathematical Society, 1997) 61--67.

\bibitem{clmno}%
F. Celler, C. R. Leedham-Green, S. H. Murray, A. C. Niemeyer and 
E.A. O'Brien, Generating random elements of a finite group, 
{\it Comm. Algebra} {\bf 23} (1995), 4931--4948.

\bibitem{cg}%
F. R. K. Chung and R. L. Graham, Random walks on generating sets of 
finite groups, {\it Electronic J. Combin.} {\bf  4(2) } (1997), 
R7 (14 pp), 
(http://www.mcs.drexel.edu /EJC/Journal/ejc-wce.html).

\bibitem{atlas}%
J.H. Conway, R.T. Curtis, S.P. Norton, R.A. Parker and 
R.A. Wilson, {\it An atlas of 
finite groups}, Clarendon Press, Oxford, 1985. 

\bibitem{demp}%
U. Dempwolff, Linear groups with large cyclic subgroups and 
translation planes, {\it Rend. 
Sem. Mat. Univ. Padova} {\bf 77} (1987), 69-113. 

\bibitem{ds1}%
P. Diaconis and L. Saloff-Coste, Walks on generating sets of abelian 
groups, {\it Probab. Theory Relat. Fields} {\bf 105} (1996), 393--421.

\bibitem{ds2}%
P. Diaconis and L. Saloff-Coste, Walks on generating sets of groups, 
{\it Technical Report, Stanford University} {\bf 497} (1996).

\bibitem{feit}%
W. Feit, On large Zsigmondy primes, {\it Proc. Amer. Math. Soc.} 
{\bf 102} (1988), 29-36.

\bibitem{gk}%
R.M. Guralnick and W.M. Kantor, Probabilistic generation of finite 
simple groups, preprint.

\bibitem{ppds}%
R. Guralnick, T. Penttila, C. E. Praeger and J. Saxl, Linear groups
with orders having certain large prime divisors, 
{\it Proc. London Math. Soc.} (to appear).

\bibitem{gs}%
R.M. Guralnick and J. Saxl, Generation of classical groups 
by conjugates, preprint.

\bibitem{her1}%
C. Hering, Transitive linear groups and linear groups which contain 
irreducible subgroups 
of prime order, {\it Geom. Ded.} {\bf 2} (1974), 425-460.  
 
\bibitem{her2}%
C. Hering, Transitive linear groups and linear groups which 
contain irreducible subgroups 
of prime order, II, {\it J. Algebra} {\bf 93} (1985), 151-164. 

\bibitem{hr}%
D. F. Holt and S. Rees, Testing modules for irreducibility,
{\it J. Austral. Math. Soc. (Series A)} {\bf 57} (1994), 1--16. 

\bibitem{james}%
G.D. James, On the minimal dimensions of irreducible 
representations of symmetric 
groups, {\it Math. Proc. Cam. Phil. Soc.} {\bf 94} (1983), 417-424. 

\bibitem{modat}%
C. Jansen, K. Lux, R. Parker and R. Wilson, 
{\it An atlas of Brauer characters}, Clarendon Press, Oxford, 1995. 

\bibitem{ks}%
E. Kaltofen and V. Shoup,
Subquadratic-time factoring of polynomials over finite fields,
{\it Math. Comp.} (to appear).

\bibitem{k-g}%
W. Keller-Gehrig,
Fast algorithms for the characteristic polynomial,
{\it Theoretical Computer Science} {\bf 36} (1985), 309--317.

\bibitem{kl}%
P. Kleidman and M. Liebeck, {\it The subgroup structure of the finite 
classical groups}, 
London Math. Soc. Lecture Note Series 129, Cambridge University Press, 
Cambridge, 1990.

\bibitem{lieb}%
M. W. Liebeck, The affine permutation groups of rank three, 
{\it Proc. London Math. Soc.} (3) {\bf 54} (1987), 477--516. 

\bibitem{lieb2}%
M. W. Liebeck, Characterization of classical groups by orbit
sizes on the natural module, preprint.


\bibitem{merkt}%
B. Merkt, {\it Zsigmondy-elemente in klassischen Gruppen}, Doctoral 
dissertation, Universit\"at zu T\"ubingen, 1995.

\bibitem{recog}%
P. M. Neumann and C. E. Praeger, A recognition algorithm for special 
linear groups, {\it Proc. London Math. Soc.} (3) {\bf 65} (1992), 555-603. 

\bibitem{dimacs}%
A. C. Niemeyer and C. E. Praeger,
Implementing a recognition algorithm for classical groups, in
{\it Groups and Computation II}, DIMACS: Series in Discrete Mathematics and
Theoretical Computer Science, {\bf 28}  (ed. L. Finkelstein and W. M.
Kantor, American Mathematical Society, 1997) 273--296.

\bibitem{recog2}%
A. C. Niemeyer and C. E. Praeger, A recognition algorithm for 
classical groups,   {\it Proc. London Math. Soc.} (to appear).
 
\bibitem{recog3}%
A. C. Niemeyer and C. E. Praeger, A recognition algorithm for 
 non-generic classical groups over finite fields, preprint, 1997.

\bibitem{nzm}%
I. Niven, H. S. Zuckerman and H. L. Montgomery,  
{\it An introduction to the  theory of  numbers}, 5th edn
(John Wiley \& Sons, New York, 1991).

\bibitem{parker}%
R. A. Parker, The computer calculation of modular
characters (the Meat-Axe), {\it Computational group theory}, Proceedings of
the London Mathematical Society Symposium on Computational Group Theory
(ed. Michael Atkinson, Academic Press, London, 1984), 267--274.

\bibitem{dim1}%
C. E. Praeger,
Computations with matrix groups over finite fields, 
{\it Groups and Computation}, DIMACS: Series in Discrete Mathematics and
Theoretical Computer Science, {\bf 11}  (ed. L. Finkelstein and W. M.
Kantor, American Mathematical Society, 1991) 189--195.

\bibitem{gap}%
M. Sch\"onert {\it et al.},  {\sf GAP}--Groups, algorithms
and programming, {\it Lehrstuhl D f\"ur Mathematik, Rheinisch Westf\"alische Technische Hochschule}, Aachen, Germany, fifth edition, 1995, 
(http://www-gap.dcs.st-and.ac.uk/~gap).

\bibitem{shalev}%
A. Shalev, A theorem on random matrices and some
applications, {\it J. Algebra} (to appear).

\bibitem{sh}%
V. Shoup, A new polynomial factorization algorithm and its implementation,
{\it J. Symbolic Comp.} {\bf 20} (1995), 363-397

\bibitem{det}%
D. E. Taylor, {\it The geometry of the classical groups}, Heldermann
Verlag, Berlin, 1992.

\bibitem{wag1}%
A. Wagner, The faithful linear representations of least degree of 
$S_n$ and $A_n$ over a 
field of characteristic 2, {\it Math. Z.} {\bf 151} (1976), 127-137. 
 
\bibitem{wag2}%
A. Wagner, The faithful linear representations of least degree of 
$S_n$ and $A_n$ over a 
field of odd characteristic, {\it Math. Z.} {\bf 154} (1977), 103-114. 

\bibitem{wag3}%
A. Wagner, An observation on the degrees of projective 
representations of the symmetric 
and alternating groups over an arbitrary field, {\it Arch. Math.} {\bf 29} (1977) 583-589.

\bibitem{zsig}%
K. Zsigmondy, Zur Theorie der Potenzreste, {\sl Monatsh. f\"ur 
Math. u. Phys.} {\bf 3} (1892), 265-284.

\end{thebibliography}
\end{document}